\documentclass[12pt,draft]{amsart}
\usepackage{amsmath,amsthm,latexsym,amscd,amsbsy,amssymb,fleqn,leqno}
\chardef\bslash=`\\ % p. 424, TeXbook

\makeatletter
\def\verbatim{\interlinepenalty\@M \@verbatim
  \leftskip\@totalleftmargin\advance\leftskip2pc
  \frenchspacing\@vobeyspaces \@xverbatim}
\makeatother
\makeatletter
  \def\dgt@k{\dg@DX=-3 \dg@DY=2 \dg@SIZE=3} 
\makeatother

\makeatletter
  \def\dgt@kk{\dg@DX=3 \dg@DY=-1 \dg@SIZE=3}%
\makeatother

\theoremstyle{plain}
\newtheorem{thm}{Theorem}[section]
\newtheorem{cor}[thm]{Corollary}
\newtheorem{lem}[thm]{Lemma}
\newtheorem{pro}[thm]{Proposition}

\newtheorem*{thrm}{Theorem}
\newtheorem*{corllr}{Corollary}

\theoremstyle{definition}

\numberwithin{equation}{section}

\newcounter{rmnum}

\newenvironment{alphanum}{\begin{list}{{\rm (\alph{rmnum})}}
{\usecounter{rmnum}\def\makelabel##1{\hss\llap{##1}}
\setlength{\leftmargin}{0pt}\setlength{\itemindent}{37pt}
\setlength{\topsep}{5pt}\setlength{\parsep}{0pt}\setlength
{\itemsep}{0pt}}}{\end{list}}

\def\symbolnote#1#2{\let\thefootn=\thefootnote%
\renewcommand{\thefootnote}{\fnsymbol{footnote}}%
\footnotemark[#1]%
\footnotetext[#1]{#2}%
\let\thefootnote=\thefootn
}

\newfont{\bbb}{msbm10 scaled \magstep1}
\newfont{\bbc}{msbm8 scaled \magstep0}

\newcommand{\N}{\mbox{\bbb N}}

\newcommand{\sphere}{\mbox{\bbb S}}
\newcommand{\uin}{\mbox{\bbb I}}
\newcommand{\e}{\mbox{\rm e-dim}}

\begin{document}

%%%%%%% Begin Topmatter %%%%%%%%%%

\title[Extension dimension and $C$-spaces]{Extension dimension and $C$-spaces}
\author{Alex Chigogidze}
\address{Department of Mathematics and Statistics,
University of Saskatche\-wan,
McLean Hall, 106 Wiggins Road, Saskatoon, SK, S7N 5E6,
Canada}
\email{chigogid@math.usask.ca}
\thanks{The first author was partially supported by NSERC grant.}

\author{Vesko Valov}
\address{Department of Mathematics, Nipissing University,
100 College Drive, P.O. Box 5200, North Bay, ON, P1B 8L7, Canada}
\email{veskov@unipissing.ca}
\thanks{The second author was partially supported by Nipissing University Research Council Grant.}

 
%%%%%%% End topmatter %%%%%%%%%

\begin{abstract}{Some generalizations of the classical Hurewicz formula are 
obtained for extension dimension and $C$-spaces. A characterization
of the class of metrizable spaces which are absolute neighborhood extensors for all
metrizable $C$-spaces is also given.}
\end{abstract}

\maketitle
\markboth{A.~Chigogidze and V.~Valov}{Extension dimension and $C$-spaces}

%%%%%%%%%%%%%%%%%%%%%%%%%%%%%%%%%%%%%%%%%%%%%%%%%%%%%%%%%%%%

\section{Introduction}

The dimension lowering Hurewicz theorem states that if
$f\colon X\to Y$ is a closed map, then 
$\dim X\leq\dim f+\dim Y$, where 
$\dim f=\sup\{\dim f^{-1}(y):y\in Y\}$ (it
was first proved by Hurewicz \cite{wh:27} for metric compacta and
later extended \cite{es:62} for paracompact spaces; 
see also \cite{bp:65}, \cite{bp:81}).
In the present paper we prove a version of Hurewicz's theorem for
extension dimension $\e$ (precise definition of this
concept is given in Section \ref{S:hur}). The ``dimesional scale"
corresponding to the extension dimension is much finer than the
usual integer-valued one. Roughly speaking extension
dimension of a space is (determined by) a complex. For instance,
the inequality
$\dim X \leq n$ is equivalent to $\e X \leq \sphere^{n}$
and the inequality $\dim_{G}X \leq n$ is equivalent to
$\e X \leq K(G,n)$ ($K(G,n)$ denotes the corresponding
Eilenberg-MacLain complex). Extension dimension allows us to detect
new properties of spaces generated by the new scale. Moreover a variety of
known facts can now be viewed from a more general point of view.

One of the first such generalizations of the classical
Hurewicz inequality was obtained in \cite{du:97}:
If $f\colon X\to Y$ is a light map (i.e. $\dim f = 0$) between
compact spaces,
then $\e X\leq\e Y$. This observation, combined
with a result of Pasynkov \cite{bp:96}, yields another
generalization of the Hurewicz formula: If $\dim f\leq n$ and $X$, $Y$
are finite-dimensional metric compacta, then 
$\e X\leq\e(Y\times\uin^n)$. The most general 
extension of the Hurewicz formula was obtained recently in \cite{drs:98}:
If $\e(Y\times f^{-1}(y))\leq K$ for every $y\in Y$, then
$\e X\leq K$ provided $X$ and $Y$ are finite-dimensional
metric compacta with $Y$ being dimensionally full-valued and
$K$ being a countable $CW$-complex. 

It is clear now that the Hurewicz formula can be generalized in 
several possible directions. One of them is to replace the
inequality $\dim f\leq n$ by $\e f^{-1}(y)\leq K$ for every
$y\in Y$ and leave $Y$ to be finite-dimensional,
say $\dim Y =m$. Note that the inequality  
$\dim X \leq n+m$ from the Hurewicz formula is equivalent to
$\e X \leq\sphere^{n+m}$
and $\sphere^{n+m}=\Sigma^m\sphere^n$, where
$\Sigma ^m\sphere^n$ denotes the $m$-iterated suspension 
of $\sphere^n$. Consequently $\dim X \leq n+m$ if and only
if $\e X \leq\Sigma^m\sphere^n$. These observations ``justify" 
our first results (for simplicity, they are not given   
in their most general forms; complete versions are recoreded
below as Theorem \ref{T1:H}
and Corollary \ref{cor2:H}) as ``extensional analogues" of
the Hurewicz formula.

\begin{thrm}\label{T1:I}
Let $f\colon X\to Y$ be a closed surjection of metrizable spaces
and $\dim Y\leq m$. 
If $K$ is a $CW$-complex such that 
$\e(\uin ^m\times f^{-1}(y))\leq K$ for any $y\in Y$, then 
$\e X\leq K$.
\end{thrm}

\begin{corllr}\label{cor1:I}
Let $f\colon X\to Y$ and the spaces $X$, $Y$ be as in the above Theorem.
Then $\e X \leq\Sigma^mK$ provided
$\e f^{-1}(y)\leq K$ for any $y\in Y$.
\end{corllr}

The second part of this paper deals with $C$-spaces \cite{ag:78}
(see also \cite{re:95}), predominantly with the class $\mathcal C$ of all
metrizable $C$-spaces.  It is well known that $\mathcal C$ contains
(strongly) countable-dimensional metrizable spaces, i.e.
metrizable spaces which are countable union of (closed)
finite-dimensional subsets,  
but there exists a metric $C$-compactum which is not
countable-dimensional \cite{p:81}. Hurewicz type theorem
is known 
\cite{hy:89} to be true for paracompact $C$-spaces (i.e.
if $f\colon X\to Y$ is a closed surjection between
paracompact spaces and if $Y$ and all fibers $f^{-1}(y)$,
$y\in Y$, are
$C$-spaces, then $X$ also is a $C$-space). Extensional
properties of $X$ in such a situation are discussed in Theorem \ref{T1:C}.
In particular, we conclude (Corollary \ref{cor1:C}) that
absolute extensors for the class $\mathcal C$, denoted by
$AE(\mathcal C)$, are precisely aspherical absolute
neighbourhood extensors for the same class ($ANE(\mathcal C)$).
Moreover in Theorem \ref{T:6} we present  description of
$ANE(\mathcal C)$-spaces and provide an answer to a
corresponding question of F. Ancel \cite[Question 5.13(c)]{fa:85}.
Another implication of Theorem 3.6 is that any subclass of
$\mathcal C$ which contains strongly countable-dimensional spaces, 
has the same absolute (neighborhood) extensors as the
class $\mathcal C$. In particular, if $\mathcal M$ is such
a proper subclass of $\mathcal C$, then $\mathcal M$ can
not be distinguished by existence of a metric space $K$
such that $\e X\leq K$ if and only if $X\in\mathcal M$
(J. Dijkstra \cite{jd:96} arrived to the same observation
for the classes $\mathcal M_\alpha$ of all metrizable spaces
with transfinite inductive dimension $\leq\alpha$, where
$\alpha$ is an infinite ordinal).

The authors are grateful to the referee whose comments and
suggestions led to a substantial improvement of the original
exposition. 
%%%%%%%%%%%%%%%%%%%%%%%%%%%%%%%%%

%%%%%%%%%%%%%%%%%%%%%%%%%%%%%%
%%%%%%%%%%%%%%%%%%%%%%%%%%%%%%

\section{Generalized Hurewicz's theorems for extension dimension}\label{S:hur}
All spaces considered in this paper are at least completely
regular and all single-valued maps are continuous. 

A space $K$ is called an absolute (neighborhood)
extensor\footnote{In these notes we follow the standard
definition of the concept of absolute (neighbourhood) extensor. It should be
noted however that in certain situations this definition is not
satisfactory and requires a modification. Such an approach is
developed in \cite{ch:00}, \cite{ch:97}, \cite{book}.} of $X$
(notation: $K\in A(N)E(X)$) if every map $f\colon A\to K$, where
defined on a closed subspace $A$ of $X$, admits an extention over
the whole $X$ (respectively, over a neighborhood of $A$ in $X$).
Next let us introduce a relation $\leq$ for $CW$-complexes. Following
\cite{d:95} (see also \cite{dd:95}, \cite{ch:00}), we say that
$L \leq K$ if for each space $X$ the condition
$L \in AE(X)$ implies the condition $K \in AE(X)$. Equivalence classes of
$CW$-complexes with respect to this relation are called extension
types. The above defined relation creates a partial order in
the collection of extension types of complexes. This partial
order is still denoted by $\leq$ and the extension type $[K]$
of a complex $K$ for simplicity is still denoted by $K$. Note
that under these definitions the collection of extension types
of all complexes has both maximal ans minimal elements. The
minimal element is the extension type of the $0$-dimensional
spehere $\sphere^{0}$ (i.e. the two-point discrete space) and
the maximal element is obviously the extension type of the
one-point space (or equivalently, of any contractible
$CW$-complex). Finally the extension dimension of a space $X$
is the minimum of extension types of complexes $K$ satisfying
the relation $K \in AE(X)$: $\e X = \min \{ [K] \colon K \in AE(X)\}$.
For simplicity below we write $\e X \leq K$ instead of
$\e X \leq [K]$.

The cone of a space $X$
(notation: $Cone(X)$) is the quotient set $X\times [0,1]/(X\times \{1\})$ with the following
topology: $U$ is open in $Cone(X)$ iff $U\cap (X\times [0,1))$ is open in $X\times [0,1)$ with the product topology and, if the vertex $v$ belongs to $U$, then $X\times (t,1)\subset U$ for some $0<t<1$. We need the following result of Dydak \cite{jd:99}: If $K$ is a space with at least two points, then $K\in ANE(X)$ if and only if $\e X \leq Cone(K)$.

\begin{lem}\label{L1:H}
Let $H\subset X$ be a zero-set in $X$ and $\e X\leq K$. Then every map 
$f\colon H\to Cone(K)\backslash\{v\}$ extends to a map from $X$ into $Cone(K)\backslash\{v\}$.
\end{lem}

\begin{proof}
Let $\pi_1\colon Cone(K)\backslash\{v\}\to K$ and $\pi_2\colon Cone(K)\to [0,1]$ be the natural projections. Then $f=(f_1,f_2)$ with $f_i=\pi_i\circ f$, $i=1,2$. Since $\e X \leq K$ implies $\e X \leq Cone(K)$ (by the Dydak result mentioned above), there exists a map $g\colon X\to Cone(K)$ extending $f$. Then $p=\pi_2\circ g$ extends $f_2$. Fix a function $q\colon X\to [0,1]$ such that $H=q^{-1}(0)$ and define $s\colon X\to [0,1)$ by $s(x)=(1-q(x))p(x)$. Since $\e X \leq K$, $f_1$ can be extended to a map $h\colon X\to K$. Then $\overline{f}=(h,s)\colon X\to Cone(K)\backslash\{v\}$ is the required extension of $f$. 
\end{proof}
Everywhere below $C(X,M)$ denotes the space of all continuous
maps from $X$ into $M$ equipped with the compact-open topology.
A set-valued map $\phi\colon X\to 2^Y$ is called strongly lower
semi-continuous (br., strongly lsc) if for  any $x\in X$
and a compact set $P\subset\phi (x)$ there exists a neighborhood
$U$ of $x$ such that $P\subset\phi (z)$ for every $z\in U$.
Here, $2^Y$ stands for the family of all nonempty subsets of $Y$.
We also write $X$ is $C^n$ to denote that every continuous image of a $k$-sphere in $X$, $k\leq n$, is contractible in $X$. 
\begin{pro}\label{P1:H}
Suppose $f\colon X\to Y$ is a closed surjection such that $X$ is
a k-space, $K\in ANE(\uin ^m\times X)$ and
$\e (\uin ^m\times f^{-1}(y))\leq K$ for any $y\in Y$.
Let $M$ be the cone of $K$ with a vertex $v$ and 
$h\colon A\to K$ a map with $A\subset X$ being a zero-set.
Then the set-valued map $\phi\colon Y\to 2^{C(X,M)}$, 
$\phi (y)=\{g\in C(X,M): g(f^{-1}(y))\subset M\backslash\{v\}
\hbox{ and $g(x)=h(x)$ for all $x\in A$}\}$ is strongly lsc
and each $\phi (y)$ is $C^{m-1}$.
\end{pro}

\begin{proof}
{\em Claim $1$.} {\em $\phi (y)\neq\emptyset$ for each $y\in Y$.}

\medskip
Observe first that $K\in ANE(\uin^m\times X)$ implies
$M\in AE(\uin^m\times X)$, in particular, $M\in AE(X)$.
For fixed $y\in Y$ extend $h|(f^{-1}(y)\cap A)$ to a map
$g_1\colon f^{-1}(y)\to K$ (such an extension exists because 
$f^{-1}(y)$ is a closed subset of $\uin^m\times f^{-1}(y))$, so $\e f^{-1}(y) \leq K$). Then
$g_1$ and $h$ define a map from $f^{-1}(y)\cup A$ into
$K$ which is extendable to a map $g\colon X\to M$.
Obviously, $g(f^{-1}(y))\subset K$ and $g|A=h|A$, so
$g\in\phi (y)$.

\medskip
{\em Claim $2$.} {\em $\phi$ is strongly lsc.}

\medskip 
Let $y_0\in Y$ and $P\subset\phi (y_0)$ be compact. We have to
find a neighborhood $V$ of $y_0$ in $Y$ such that 
$P\subset\phi (y)$ for every $y\in V$. Let
$P(x)=\{g(x):g\in P\}$, $x\in X$. Since $P\subset C(X,M)$ is compact 
and $X$ is a $k$-space, by the Ascoli theorem, each $P(x)$ 
is compact and $P$ is evently continuous. This easily implies
that the set
$W=\{x\in X: P(x)\subset M\backslash\{v\}\}$
is open in $X$ and, obviously, $f^{-1}(y_0)\subset W$. 
Because $f$ is closed, there exists a neighborhood $V$
of $y_0$ in $Y$ with $f^{-1}(V)\subset W$. Then, according
to the choice of $W$ and the definition of $\phi$, 
$P\subset\phi (y)$ for every $y\in V$.

\medskip
{\em Claim $3$.} {\em Each $\phi (y)$ is $C^{m-1}$}.

\medskip
For a fixed $y\in Y$ take an arbitrary map
$u\colon\sphere ^{n-1}\to\phi (y)$, where $n\leq m$.
We are going to show that $u$ can be
extended continuously to a map from $\uin ^n$ into $\phi (y)$
(we identify $\sphere ^{n-1}$ with the boundary of $\uin ^n$).
Since 
$\sphere ^{n-1}\times X$ is a k-space (as a product of a compact space
and a k-space), the map $u_1\colon\sphere ^{n-1}\times X\to M$, 
$u_1(z,x)=u(z)(x)$, is continuous (see \cite{re:89}).
Because $u_1(z,x)=h(x)$ for every $(z,x)\in\sphere ^{n-1}\times A$,
we can extend $u_1|(\sphere ^{n-1}\times A)$ to a map 
$u_2\colon\uin ^n\times A\to K$, $u_2(z,x)=h(x)$. 
Then, we have a closed subset
$H=(\sphere^{n-1}\times f^{-1}(y))\cup (\uin ^n\times (f^{-1}(y)\cap A))$
of $\uin ^n\times f^{-1}(y)$ and a map 
$u_3\colon H\to M\backslash\{v\}$ defined by
$u_3|(\sphere^{n-1}\times f^{-1}(y))=u_1|(\sphere^{n-1}\times f^{-1}(y))$ and
$u_3|(\uin ^n\times (f^{-1}(y)\cap A))=u_2|(\uin ^n\times (f^{-1}(y)\cap A))$.
Since $\sphere^{n-1}$ and $f^{-1}(y)\cap A$ are zero-sets in $\uin ^n$ and $f^{-1}(y)$,
respectively, both $\sphere^{n-1}\times f^{-1}(y)$ and $\uin^n\times (f^{-1}(y)\cap A)$ are zero-sets in $\uin^n\times f^{-1}(y)$, so is $H$.
Note that $\e(\uin^n\times f^{-1}(y))\leq K$ because $\uin^n\times f^{-1}(y)$ is closed in $\uin^m\times f^{-1}(y)$.
Therefore, by Lemma 2.1, 
$u_3$ extends to a map
$u_4\colon\uin ^n\times f^{-1}(y)\to M\backslash\{v\}$. Now, let
$F$ be the union of the sets $F_1=\uin ^n\times f^{-1}(y)$,
$F_2=\uin ^n\times A$ and $F_3=\sphere ^{n-1}\times X$.  We define
the map $p\colon F\to M$ by $p|F_1=u_4$, $p|F_2=u_2$ and $p|F_3=u_1$.
Obviously, $F$ is closed in $\uin ^n\times X$.
Since $M\in AE(\uin ^n\times X)$, there exists an extension
$q\colon\uin ^n\times X\to M$ of $p$. To finish the proof of Claim 3,
observe that $q$ generates the map $\overline{u}\colon\uin ^n\to C(X,M)$,
$\overline{u}(z)(x)=q(z,x)$. Moreover, $q(z,x)=h(x)$ for any
$(z,x)\in\uin ^n\times A$ and
$q(\uin ^n\times f^{-1}(y))\subset M\backslash\{v\}$.
So, $\overline{u}$ is a map from $\uin ^n$ to $\phi (y)$ which extends $u$.
\end{proof}

Now we need the following result of E. Michael \cite[Remark 2]{vu:98}.

\begin{pro}\label{M}
Let $X$ be paracompact with $\dim X\leq m$ and $Y$ an arbitrary space.
Then every strongly lsc mapping $\varphi\colon X \to 2^Y$    
has a continuous selection provided $\varphi (x)$ is $C^{m-1}$
for each $x\in X$.
\end{pro}

\begin{thm}\label{T1:H}
Let $f\colon X\to Y$ be a closed surjection with $X$ a k-space and
$Y$ paracompact of dimension $\dim Y\leq m$.
If $K$ is any space such that $K\in ANE(\uin ^m\times X)$ and
$\e (\uin ^m\times f^{-1}(y))\leq K$ for any $y\in Y$, then $\e X \leq K$.
\end{thm}

\begin{proof}
Suppose $A\subset X$ is closed and $h\colon A\to K$ is a map.
We are going to find a continuous extension $\overline{h}\colon X\to K$
of $h$. Let $M$ be the cone of $K$ with a vertex $v$. Since $M\in AE(X)$, there exists a map
$q\colon X\to M$ extending $h$. Then $q^{-1}(K)$ is a zero-set in $X$ (because $K$ is such a set in $M$) containing $A$. Therefore, we can assume that $A$ is a zero-set in $X$. Next,
define the set-valued map
$\phi\colon Y\to 2^{C(X,M)}$,
$\phi (y)=\{g\in C(X,M): g(f^{-1}(y))\subset M\backslash\{v\}
\hbox{ and $g(x)=h(x)$ for all $x\in A$}\}$ (a similar idea was
earlier used by V.Gutev and V. Valov).
By Proposition \ref{P1:H},
$\phi\colon Y\to C(X,M)$ is a strongly lsc map with each $\phi (y)$
being a $C^{m-1}$-set. Since $\dim Y\leq m$, we can apply
Proposition \ref{M}
to obtain a continuous selection $t\colon Y\to C(X,M)$ for $\phi$.
Then $g\colon X\to M$, defined by $g(x)=t(f(x))(x)$, is 
continuous on every compact subset of $X$ and because
$X$ is a k-space, $g$ is continuous. Since $t(f(x))\in\phi (f(x))$,
we have $g(x)=h(x)$ for all $x\in A$ and $g(x)\in M\backslash\{v\}$,
$x\in X$. Finally, if $\pi_1\colon M\backslash\{v\}\to K$
denotes the natural retraction, then
$\overline{h}=\pi\circ g\colon X\to K$ is the required continuous
extension of $h$.
\end{proof}

A $k$-space $X$ is called a $cw$-space \cite{dd:95} if every contractible $CW$-complex is an $AE(X)$. In particular, if $X$ is a $cw$-space and $K$ any $CW$-complex, then $Cone(K)\in AE(X)$. Any metrizable space, more generally, every space admitting a perfect map onto a first countable paracompact space, is $cw$ \cite{jd:99}.
\begin{cor}\label{cor1:H}
Let $f\colon X\to Y$ be a closed surjection, where  
$Y$ is paracompact with $\dim Y \leq m$ and $\uin^m\times X$ is
a cw-space. If $K$ is a $CW$-complex such that 
$\e(\uin ^m\times f^{-1}(y))\leq K$ for every $y\in Y$, then $\e X \leq K$.
\end{cor}

\begin{proof}
Since $X$ is a k-space and $K\in ANE(\uin^m\times X)$,
we can apply Theorem \ref{T1:H}.
\end{proof}

\begin{lem}\label{L2:H}
If $\e X \leq K$, where $X\times\uin$ is a paracompact cw-space and $K$
a $CW$-complex, then $\e(X\times\uin)\leq\Sigma K$.
\end{lem}

\begin{proof} This lemma was proved by Dranishnikov \cite{ad:91} for metric
spaces $X$. His proof, coupled with \cite[Propositions 1.17-1.18]{dd:95}, 
works in our situation as well.
\end{proof}

\begin{cor}\label{cor2:H}
Let $X\times\uin^m$ be a paracompact cw-space, $K$ be a $CW$-complex and
$f\colon X\to Y$ be a closed surjection with $\dim Y\leq m$.
If $\e f^{-1}(y)\leq K$ for every $y\in Y$, then $\e X \leq\Sigma ^mK$. 
\end{cor}

\begin{proof}
Observe first that $Y$ is paracompact as a closed image of the paracompact $X$.
By Lemma \ref{L2:H}, $\e(\uin ^m\times f^{-1}(y))\leq\Sigma ^mK$
for any $y\in Y$. Then the proof follows from Corollary \ref{cor1:H} with $K$
replaced by $\Sigma ^mK$.
\end{proof}. 

\section{$C$-spaces}

Recall that $X$ is a $C$-space \cite{ag:78} if
for any sequence $\{\omega _n\}$ of open covers of $X$ there
exists a sequence
$\{\gamma _n\}$ of open disjoint families in $X$ such that each $\gamma _n$
refines $\omega _n$ and $\bigcup\{\gamma _n:n\in\N\}$ covers $X$. 
Property $C$ is a dimensional type
property, and it admits a characterization
similar to that one (see Proposition \ref{M}) of finite-dimensional spaces
(everywhere below a space is said to be aspherical if it is $C^n$ for all $n$).

\begin{pro}$[20]$\label{U}
A paracompact $X$ is a $C$-space if and only if every strongly
lsc map $\phi\colon X\to 2^Y$ with aspherical images $\phi (x)$, $x\in X$,
where $Y$ is an arbitrary space,
has a continuous selection.
\end{pro}

\begin{thm}\label{T1:C}
Let $f\colon X\to Y$ be a closed surjection with $X$ a k-space and
$Y$ a paracompact $C$-space.
If $K$ is a space satisfying both conditions $K\in ANE(\uin ^m\times X)$  and
$K \in AE(\uin ^m\times f^{-1}(y))$ for any $m\in\N$ and any $y\in Y$,
then $K \in AE(X)$.
\end{thm}

\begin{proof}
We follow the proof of Theorem \ref{T1:H}. Maintaining the same
notations and applying now Proposition \ref{U} (instead of
Proposition \ref{M}), it suffices to show that if $A$ is a zero-set in $X$, then
the formula   
$\phi (y)=\{g\in C(X,M): g(f^{-1}(y))\subset M\backslash\{v\}\hbox{
and $g(x)=h(x)$
for all $x\in A$}\}$ defines a set-valued map
$\phi\colon Y\to 2^{C(X,M)}$ which is strongly lsc
and each $\phi (y)$ is aspherical. And this follows from Proposition \ref{P1:H}.
\end{proof}
Theorem 3.2 is not of any interest when $K$ is a $CW$-complex. Indeed,  
$K \in AE(\uin^m\times f^{-1}(y))$ for all $m$ implies that every homotopy group of $K$ is trivial. So, $K$ is contractible and therefore it is an absolute extensor for any $cw$-space.
On the other hand, the Borsuk example of a contractible and locally contractible compact metric space which is not an $AE$ for the class of all metrizable spaces shows that Theorem 3.2 has a meaning for general spaces $K$.

Let $\mathcal C$ denote the class of all metrizable $C$-spaces. We write $K\in A(N)E({\mathcal C})$ if $K\in A(N)E(X)$ for any $X\in{\mathcal C}$; when the class of all metrizable spaces is considered, we simply write $K\in A(N)E$. 

\begin{cor}\label{cor1:C}
A space $K\in ANE({\mathcal C})$ is an $AE({\mathcal C})$ if and only if $K$ is aspherical.
\end{cor}

\begin{proof}
Any $AE(\mathcal C)$ is aspherical (because the class $\mathcal C$ contains all finite-dimensional spaces) and an $ANE(\mathcal C)$.
Suppose $K\in ANE({\mathcal C})$ is aspherical and $X\in{\mathcal C}$. We are going to apply
Theorem 3.2 in the special case when $X=Y$ and $f$ being the identity map. In this special case Proposition 2.2 is true if $\e(\uin^m\times f^{-1}(y))\leq K$ is replaced by 
$K\in C^{m-1}$. Indeed, Claim 1 becomes trivial; to prove Claim 2 we don't need to apply Lemma 2.1 because the set $H$ is homeomorphic either to $\uin^n$ if $y\in A$ or $\sphere^{n-1}$ otherwise, we need that any map from $\sphere^{n-1}$ into $K$ is extendable to map 
from $\uin^n$ into $K$, $n\leq m$. In order to apply Theorem 3.2, it remains only to check that
$K\in ANE(X\times\uin^m)$ for all $m$. And that is true because $X\times\uin^m\in\mathcal C$ \cite{hy:89}. 
\end{proof}

Let discuss now some sufficient (and necessary) conditions for a metric space to be an $ANE(\mathcal C)$. Let $\mathcal P$ be a topological property. We say that
$X\subset E$ is a $UV(\mathcal P)$ subset of $E$ if each neighborhood $U$ of $X$ in $E$ contains a neighborhood $V$ of $X$ in $E$ such that any map $h\colon Z\to V$, where $Z\in\mathcal P$, extends to a map $\overline{h}\colon Cone(Z)\to U$. 
A closed surjection $f\colon X\to Y$ is called $UV(\mathcal P)$ if each of its point inverses is a $UV(\mathcal P)$ subset of $X$. Recall that if, in the above definition, $V$ is contractible in $U$, then $X$ is called $UV^{\infty}$; a cell-like space is a compact metric space $X$ such that $X$ is a $UV^{\infty}$ set in every $ANE$-space $E$ in which it is embedded as a closed subset (see, for example, \cite{fa:85}). In the existing terminology, a $UV^{\infty}$ (resp., cell-like) map is a perfect map with $UV^{\infty}$ (cell-like) preimages. Obviously, every $UV^{\infty}$ map is $UV(\mathcal P)$ for any property $\mathcal P$.     

\begin{pro}
Any one of the following two conditions is sufficient for a metrizable
space $Y$ to be an $ANE(\mathcal C)$:
\begin{alphanum}
\item
$Y$ is locally contractible, more generally, there exists a metrizable space $X$ 
and a $UV^{\infty}$ map from $X$ onto $Y$. 
\item
$Y$ has a base of open aspherical sets.
\end{alphanum}
\end{pro}

\begin{proof}
First condition was proved by Ancel \cite[Theorem C.5.9]{fa:85}, see also \cite{ag:78} for the case of local contractibility. Condition (b) can be obtained by using the arguments of Ageev and Repov\v{s}
\cite[proof of Theorem 1.3]{ar:99}.
\end{proof}

Not every metrizable $ANE(\mathcal C)$-space is locally contractible. J. van Mill provided an example of a cell-like image of the Hilbert cube such that no nonempty open subset is contractible in that space \cite{vm:86}. At the same time, by \cite{fa:85}, this example is an $ANE(\mathcal C)$. In view of mentioned above result of Ancel \cite[Theorem C.5.9]{fa:85}, it is interesting whether any metrizable $ANE(\mathcal C)$ is a $UV^{\infty}$ image of a metrizable space. In such a case, the class of metrizable $ANE(\mathcal C)$ would be precisely the class of all $UV^{\infty}$ images of metrizable spaces.  
We can provide similar characterization $ANE(\mathcal C)$ in terms of $UV(s.c.d.)$ maps, where     
s.c.d. denotes the property strong countable-dimensionality.

\begin{pro}
Let $f\colon M\to X$ be a surjective map between metrizable spaces. If for any $x\in X$ and its neighborhood $U(x)$ in $X$ there exists another neighborhood $V(x)$ of $x$ in $X$ such that $\overline{V}(x)=f^{-1}(V(x))$ is contractible in $\overline{U}(x)=f^{-1}(U(x))$, then $X\in ANE(\mathcal C)$.
\end{pro}

\begin{proof}
First step is to show that $X$ is an approximate absolute neighborhood extensor for the class $\mathcal C$, i.e. if $H$ is a metrizable $C$-space, $A\subset H$ closed and $h\colon A\to X$ a map, then for every open cover $\gamma$ of $X$ there is a neighborhood $W_A$ of $A$ in $H$ and a map $\overline{h}\colon W_A\to X$ such that
$\overline{h}|A$ is $\gamma$-close to $h$. We follow the construction from the proof of  \cite[Teorem 4.3, first part]{ar:99}. For every $x\in X$ and $n\geq 0$ fix points $z(x)\in f^{-1}(x)$ and neighborhoods $V_n(x)\subset U_n(x)$ of $x$ in $X$ such that:

\smallskip
\begin{alphanum}
\item[(1)]
$\overline{V}_n(x)$ contracts in $\overline{U}_n(x)$ to $z(x)$ 
for all $n\geq 0$ and $x\in X$;

\item[(2)]
the cover $\alpha_0=\{U_0(x):x\in X\}$ refines $\gamma$;
\item[(3)]
the cover $\alpha_n=\{U_n(x):x\in X\}$ star-refines $\beta_{n-1}=\{V_{n-1}(x):x\in X\}$ for any $n\geq 1$, i.e. $\{St(U,\alpha_n):U\in\alpha_n\}$ refines $\beta_{n-1}$.
\end{alphanum}

Observe that we have corresponding covers $\overline{\gamma}=f^{-1}(\gamma)$, 
$\overline{\alpha}_n=\{\overline{U}_n(x):x\in X\}$ and 
$\overline{\beta}_n=\{\overline{V}_n(x):x\in X\}$ of $M$ such that 
$\overline{\alpha}_0$ refines $\overline{\gamma}$ and
$\overline{\alpha}_n$ star-refines $\overline{\beta}_{n-1}$, $n\geq 1$. 
For every $n\geq 0$ and $x\in X$ we fix a contraction map 
$F^{x,n}\colon\overline{V}_n(x)\times [0,1]\to\overline{U}_n(x)$ with $F^{x,n}(z,1)=z(x)$. 
Since $A$ is a $C$-space (as a closed subset of $H$), there is a sequence of disjoint open families $\{\mu_n:n=1,2,..\}$ in $H$ such that the restriction of each $\mu_n$ on $A$ refines $h^{-1}(\beta_n)$ and $\mu =\bigcup\{\mu_n:n=1,2..\}$ covers $A$. Further, let $\mathcal K$ be the nerve of $\mu$ and $\theta\colon W_A=\cup\{W:W\in\mu\}\to |\mathcal K|$ a barycentric map. 
We are going to define a map $g\colon |{\mathcal K}|\to M$ such that the family 
$\{g(\theta (y))\cup f^{-1}(h(y)):y\in A\}$ refines $\overline{\gamma}$. Then the map
$\overline{h}=f\circ g\circ\theta$ will be the required $\gamma$-approximation of $h$.
Any simplex $(W_0,W_1,..,W_k)$ from $\mathcal K$, where $W_i\in\mu_{n(i)}$, can be ordered
such that $n(0)<n(1)<...,n(k)$ (this is possible because $\cap\{W_i:i=1,2,..,k\}\neq\emptyset$, so the numbers $n(i)$ are different). By $(3)$,
for any $W\in\mu_n$ there exists $x(W)\in X$ with $St(h(W\cap A),\alpha_n)\subset V_{n-1}(x(W))$. We define $g_0\colon |{\mathcal K}^0|\to M$ by $g_0(W)=z(x(W))$, $W\in\mu$. Using the contractions $F^{x,n}$, as in \cite[proof of Theorem 4.3]{ar:99}, we can define by induction maps $g_n\colon |{\mathcal K}^n|\to M$ such that the restriction of $g_n$ on $|{\mathcal K}^i|$ is $g_i$, $i\leq n$, and for any simplex $\bigtriangleup ^n=(W_0,W_1,..,W_n)\in |{\mathcal K}^n|$ we have

\smallskip
\begin{alphanum}
\item[(4)]
$f^{-1}(h(W_0\cap A))\cup g_n(\bigtriangleup ^n)\subset\overline{U}_{n_0-1}(x(W_0))$.
\end{alphanum}

\smallskip  
So, we obtain a map $g\colon |{\mathcal K}|\to M$ and, by $(4)$, $\overline{h}|A$ and $h$ are $\gamma$-close, where $\overline{h}=f\circ g\circ\theta$. Indeed, if $y\in A$ and $\theta (y)\in\bigtriangleup ^n$ for some simplex $\bigtriangleup ^n=(W_0,W_1,..,W_n)$, then
$\overline{h}(y)\in f(g_n(\bigtriangleup ^n))$ and $h(y)\in h(W_0\cap A)$. According to $(4)$, the last two inclusions imply that both $\overline{h}(y)$ and $h(y)$ belong to 
$U_{n_0-1}(x(W_0))$. So, $\overline{h}(y)$ and $h(y)$ are $\alpha_{n_0-1}$-close and, since $n_0-1\geq 0$, they are also $\gamma$-close. Therefore, $X$ is an approximate absolute neighborhood extensor for the class $\mathcal C$. 

To complete the proof we state the following result which was actually proved in \cite{ar:99} but not explicitely formulated: If $\mathcal M$ is a class of metrizable spaces such that $Y\times [0,1)\in\mathcal M$ for every $Y\in\mathcal M$, then any approximate absolute neighborhood extensor for $\mathcal M$ is an $ANE(\mathcal M)$. Since $\mathcal C$ is closed with respect to multiplication by $[0,1)$, we have $X\in ANE(\mathcal C)$.       
\end{proof} 

\begin{thm}\label{T:6}
For a metrizable space $X$ the following conditions are equivalent:
\begin{alphanum}
\item
$X$ is an $ANE$ for the class of metrizable (strongly) countable-dimensional spaces.
\item
$X$ is a $UV(s.c.d.)$ image of a metrizable space. 
\item
$X$ is an $ANE(\mathcal C)$.
\end{alphanum}
\end{thm}

\begin{proof}
Since every metrizable (strongly) countable-dimensional space has property $C$, (c) implies (a). Standard arguments show that every metrizable $X$ which is an
$ANE$ for the class of metrizable (strongly) countable-dimensional spaces has the following property $(*)$:

For every $x\in X$ and its neighborhood $U(x)$ in $X$ there is a neighborhood $V(x)\subset U(x)$ such that any map from a closed subset of a (strongly) countable-dimensional metrizable space $Z$ into $V(x)$ extends to a map from $Z$ into $U(x)$. 

Hence, $(a)$ yields that the identity map of $X$ is $UV(s.c.d)$.
So, it remains to prove (b)$\Rightarrow$(c). Let $f\colon Y\to X$ be a $UV(s.c.d.)$ map with $Y$ metrizable.
We need the following result of M. Zarichnyi \cite{mz:95}: There exists an $\omega$-soft map from a $\sigma$-compact strongly countable-dimensional metrizable space onto the Hilbert cube. Here, a map $g\colon M\to H$ is called $\omega$-soft if for every strongly countable-dimensional metrizable space $Z$, its closed subset $B\subset Z$ and any two maps $\phi\colon Z\to H$, $\psi\colon B\to M$ such that $g\circ\psi =\phi |B$ there exists a map $\Phi\colon Z\to M$ extending $\psi$ with $g\circ\Phi =\phi$. Using the Zarichnyi result, for every cardinal $\tau$ we can construct a strongly countable-dimensional metrizable space $M(\tau)$ of weight $\tau$ and an $\omega$-soft map 
$g\colon M(\tau)\to l_2(\tau)$ (see \cite{cv:90} for a similar reduction), where $l_2(\tau)$ denotes the Hilbert space of weight $\tau$. Embedding $Y$ into $l_2(\tau)$  for some $\tau$ and considering the restriction $g_Y$ of $g$ onto $M_Y=g^{-1}(Y)$, we obtain a strongly countable-dimensional metrizable space $M_Y$ and an $\omega$-soft map $g_Y\colon M_Y\to Y$.  
Let $q=f\circ g_Y$. We are going to show that $q\colon M_Y\to X$ satisfies the hypotheses of Proposition 3.5. To this end, let $U(x)$ be a neighborhood of $x\in X$. Since $f$ is $UV(s.c.d.)$, there exists a neighborhood $W(x)\subset f^{-1}(U(x))$ such that every map from a strongly countable-dimensional metrizable space $Z$ into $W(x)$ extends to a map from $Cone(Z)$ into $f^{-1}(U(x))$. Then $f^{-1}(V(x))\subset W(x)$ for some neighborhood $V(x)$ of $x$ in $X$ because $f$ is closed. Now consider $\overline{V}(x)=q^{-1}(V(x))$ and 
$\overline{U}(x)=q^{-1}(U(x))$. Since $\overline{V}(x)$ is strongly countable-dimensional, there exists a map $\phi\colon Cone(\overline{V}(x))\to f^{-1}(U(x))$ extending the restriction $g_Y|\overline{V}(x)$. Finally, using that $g_Y$ $\omega$-soft, we can lift $\phi$ to a map $\Phi\colon Cone(\overline{V}(x))\to\overline{U}(x)$ such that 
$\Phi |\overline{V}(x)$ is the identity. Therefore, $\overline{V}(x)$ is contractible in
$\overline{U}(x)$ and, by Proposition 3.5, $X\in ANE(\mathcal C)$.      
\end{proof}

The equivalence of conditions $(a)$ and $(c)$ from Theorem 3.6, yields the following observation: if $\mathcal M$ is a subclass of $\mathcal C$ containing all strongly countable-dimensional spaces, then $ANE(\mathcal M)$ coincides with $ANE(\mathcal C)$ in the realm of metrizable spaces. Consequently, since every $K\in AE(\mathcal M)$ is aspherical, the above observation combined with Corollary 3.3 implies also that $\mathcal M$ and $\mathcal C$ have the same metrizable $AE$-spaces. Finally, we would like to point out that Theorem 3.6 provides an answer to the question \cite[Question 5.13(c)]{fa:85} asking whether a metrizable space $X$ is an $ANE$ for the class of countable-dimensional spaces if $X$ has the property $(*)$ mentioned in the proof of Theorem 3.6.

\bigskip


\begin{thebibliography}{99}

\bibitem{ag:78}
D.~Addis and J.~Gresham, {\em A class of infinite-dimensional spaces. Part I:
Dimension theory and Alexandroff's Problem}, Fund. Math.
{\bf 101} (1978), 195--205.

\bibitem{ar:99}
S.~Ageev and D.~Repov\v{s}, {\em A method of approximate extension of maps in theory of extensors}, University of Ljubljana, Preprint series, Vol. {\bf 37} (1999), 651.

\bibitem{fa:85}
F.~Ancel, {\em The role of countable dimensionality in the theory of cell-like relations},
Trans. Amer. Math. Soc. {\bf 287} (1985), 1--40.

\bibitem{ch:00}
A.~Chigogidze, {\em Infinite dimensional topology and shape theory}, in:
Handbook of Geometric Topology, ed. by R. Daverman and R. Sher, North Holland,
Amsterdam, 2000 (to appear).

\bibitem{ch:97}
\bysame, {\em Cohomological dimension of Tychonov spaces}, Topology and Appl.
{\bf 79} (1997), 197--228.

\bibitem{book}
\bysame, {\em Inverse Spectra}, North Holland, Amsterdam, 1996.

\bibitem{cv:90}
\bysame and V.~Valov {\em Universal maps and surjective characterizations of completely metrizable $LC^n$-spaces}, Proc. Amer. Math. Soc. {\bf 109} (1990), 1125--1133.

\bibitem{jd:96}
J.~Dijkstra, {\em A dimension raising hereditary shape equivalence}, Fund. Math. {\bf 149}
(1996), 265--274.

\bibitem{ad:91}
A.~Dranishnikov, {\em On intersection of compacta in Euclidean space II},
Proc. Amer. Math. Soc. {\bf 113, 4} (1991), 1149--1154.

\bibitem{d:95}
\bysame, {\em The Eilenberg-Borsuk theorem for mappings into an
arbitrary complex}, Russian Acad. Sci. Sb. {\bf 81} (1995), 467--475.

\bibitem{dd:95}
\bysame and J.~Dydak, {\em Extension dimension and extension types}, Proceedings of the 
Steklov Inst. of Math. {\bf 212} (1996), 55--88.



\bibitem{drs:98}
\bysame, D.~Repov\v{s} and E.~\v{S}\v{c}epin,
{\em Transversal intersection formula for compacta}, Topology and Appl.
{\bf85} (1998), 93--117.

\bibitem{du:97}
\bysame and V.~Uspenskij, {\em Light maps and extensional dimension}, Topology and Appl. 
{\bf 80} (1997), 91--99.

\bibitem{jd:99}
J.~Dydak, {\em Extension theory: The interface between set-theoretic and
algebraic topology}, Topology Appl. {\bf 20} (1996), 1--34.

\bibitem{re:89}
R.~Engelking, {\em General Topology} (Heldermann Verlag, Berlin, 1989).

\bibitem{re:95}
\bysame, {\em Theory of dimensions: Finite and Infinite}
(Heldermann Verlag, Lemgo, 1995).

\bibitem{hy:89}
Y.~Hattori and K.~Yamada, {\em Closed preimages of $C$-spaces},
Math. Japon. {\bf 34} (1989), 555--561.

\bibitem{wh:27}
W.~Hurewicz, {\em Uber Stetige Bilder von Punktmengen (Zweite
Mittelung)}, Proc. Akad. Amsterdam {\bf 30, 1} (1927), 159--165.

\bibitem{vm:86}
Jan van Mill, {\em Local contractibility, cell-like maps, and dimension}, Proc. Amer. Math. Soc. {\bf 98} (1986), 534--536.
 
\bibitem{bp:65}
B.~Pasynkov, {\em On the Hurewicz formula}, Vestnik Mosk. Univ.
Ser. Mat. {\bf 4} (1965), 3--5 (in Russian).

\bibitem{bp:81}
\bysame, {\em Factorization theorems in dimension theory}, Uspehi Mat. Nauk
{\bf 36, 3} (1981), 147--175 (in Russian).

\bibitem{bp:96}
\bysame, {\em On geometry of continuous maps of finite-dimensional metric
compacta}, Trudy Steklov Math. Inst. {\bf 212} (1996),
147--172 (in Russian).

\bibitem{p:81}
R.~Pol, {\em A weakly infinite-dimensional compactum which is not
countable-dimensional}, Proc. Amer. Math. Soc. {\bf 82} (1981), 634--636.

\bibitem{es:62}
E.~Skljarenko, {\em A theorem on dimension-lowering mappings},
Bull. Acad. Pol. Sci. Ser. Math. {\bf 10} (1962), 429--432 (in Russian).

\bibitem{vu:98}
V.~Uspenskij, {\em A selection theorem for $C$-spaces}, Topology and Appl. {\bf 85}
(1998), 351--374.

\bibitem{mz:95}
M.~Zarichnyi, {\em Universal map of $\sigma$ onto $\Sigma$ and absorbing sets in the classes of absolute Borelian and projective finite-dimensional spaces}, Topology and Appl. {\bf 67} (1995), 221--230.

\end{thebibliography}
\end{document}